\documentclass{amsproc}
\usepackage{amsthm,amsmath,amsfonts,mathrsfs,amssymb}
\title{The centralizer of an $I$-matrix in $M_2(R/I)$,$\;R$ a UFD}
\author{Magdaleen S. Marais}
\email{magdaleen@aims.ac.za}
\address{African Institute for Mathematical Sciences, 6 Melrose Rd, Muizenberg, 7945, Cape Town, South Africa }
\thanks{ This research is part of the author's research for her doctoral dissertation which was conducted at Stellenbosch University under the direction of L. van Wyk. The financial assistance of the National Research Foundation (NRF) towards this research is hereby acknowledged. Opinions expressed and conclusions arrived at are those of the author and are not necessarily to be attributed to the National Research Foundation.}
\begin{document}
\begin{abstract}
The concept of an $I$-matrix in the full $2\times 2$ matrix ring $M_2(R/I)$, where $R$ is an arbitrary UFD and $I$ is a nonzero ideal in~$R$, is introduced. We obtain a concrete description of the centralizer of an $I$-matrix~$\widehat B$ in $M_2(R/I)$ as the sum of two subrings $\mathcal S_1$ and $\mathcal S_2$ of $M_2(R/I)$, where~$\mathcal S_1$ is the image (under the natural epimorphism from $M_2(R)$ to $M_2(R/I)$) of the centralizer in $M_2(R)$ of a pre-image of $\widehat B$, and where the entries in $\mathcal S_2$ are intersections of certain annihilators of elements arising from the entries of~$\widehat B$. It turns out that if $R$ is a PID, then every matrix in $M_2(R/I)$ is an $I$-matrix. However, this is not the case if $R$ is a UFD in general. Moreover, for every factor ring $R/I$ with zero divisors and every $n\ge 3$ there is a matrix for which the mentioned concrete description is not valid.
\end{abstract}
\maketitle
\section{Introduction}
\noindent We denote the centralizer of an element $s$ in an arbitrary ring $S$ by Cen$_{S}(s)$. Knowing that $M_n(R)$, the full $n\times n$ matrix ring over a commutative ring $R$, is a prime example of a non-commutative ring, it is surprising that a concrete description of $\textnormal{Cen}_{M_n(R)}(B)$ for an arbitrary $B\in M_n(R)$ has not yet been found. If $R[x]$ is the polynomial ring in the variable $x$ over $R$, then
\begin{equation}\label{1a}
\{f(B)\ |\ f(x)\in R[x]\}\subseteq \textrm{Cen}_{M_n(R)}(B).
\end{equation}
In fact, it is known that (see \cite{cen})
\begin{equation*}
\{f(B)\ |\ f(x)\in R[x]\}=\textrm{Cen}_{M_n(R)}(\textrm{Cen}_{M_n(R)}(B)).
\end{equation*}
The most progress, finding a concrete description of $\textnormal{Cen}_{M_n(R)}(B)$, has been made for the case when the underlying ring $R$ is a field (see \cite{halm}, \cite{hung}, \cite{frob}, \cite{closed bases} and \cite{sup}).
The following well-known result in this case provides a necessary and sufficient condition for equality in (\ref{1a}).

\newtheorem{Theorem2.1}{\bf Theorem}[section]
\begin{Theorem2.1}\label{Theorem2.1}
If $B$ is an $n\times n$ matrix over a field $F$, then
\[\textnormal{Cen}_{M_n(F)}(B)=\{f(B)\ |\ f(x)\in F[x]\}\]
if and only if the minimum polynomial of $B$ coincides with the characteristic polynomial of
$B$.
\end{Theorem2.1}

In this paper we consider the centralizer of a so-called $I$-matrix in $M_2(R/I)$, with $R/I$ a factor ring of a UFD $R$ and $I$ a nonzero ideal in $R$.

In Section~2 we obtain an explicit description of the centralizer of a $2\times 2$ matrix over a field or over a unique factorization domain. Section~2 also contains other preliminary results concerning the centralizer of an $n\times n$ matrix that will be used in the subsequent sections, including Proposition \ref{Lemma5} which may be considered as the inspiration behind this paper. In this proposition we show that the centralizer of an~$n\times n$ matrix~$\widehat B$ over a homomorphic image $S$ of a commutative ring $R$ contains the sum of two subrings $\mathcal S_1$ and $\mathcal S_2$ of $M_2(S)$, where $\mathcal S_1$ is the image of the centralizer in~$M_2(R)$ of a pre-image of $\widehat B$, and where the entries in $\mathcal S_2$ are intersections of certain annihilators of elements arising from the entries of $\widehat B$. 

In Section~3 we introduce the concepts of $I$-invertibility in a factor ring $R/I$ of a UFD $R$ (Definition~\ref{Def8}) and of an $I$-matrix in $M_2(R/I)$ (Definition~\ref{Def10a}). We show in Corollaries \ref{Corollary1} and \ref{Corollary2} that if $R$ is a PID, then every element in~$R/I$ is $I$-invertible and every matrix in $M_2(R/I)$ is an $I$-matrix. Examples \ref{Example2b} and~\ref{Example}(b) show that this is not true for UFD's in general, not even if $I$ is a principal ideal.

Section~4 contains the main result of the paper, namely Theorem \ref{Theorem21}, which provides a concrete description of the centralizer of an $I$-matrix in $M_2(R/I)$ as the sum of the above mentioned two subrings, where $R$ is a UFD and $I$ is a nonzero ideal in~$R$.

Since every $2\times 2$ matrix over a factor ring of a PID is an $I$-matrix, Theorem~\ref{Theorem21} applies to all $2\times 2$ matrices over factor rings of PID's. In Example~\ref{Example20} we exhibit a UFD~$R$, which is not a PID, a finitely generated ideal $I$ and a matrix in $M_2(R)$, which is not an~$I$-matrix, for which Theorem~\ref{Theorem21} does not hold. In Example~\ref{Example21b} we show that if $R$ is a UFD and $R/I$ is such that $R/I$ is not an integral domain, then for every~$n\ge 3$ there is a matrix in $M_n(R)$ for which we do not have equality in Proposition~\ref{Lemma5}.

\section{Preliminary Results}
\noindent Since the minimum polynomial and characteristic polynomial of any $2\times 2$ non-scalar matrix over a field coincide, the following corollary follows from Theorem~\ref{Theorem2.1}:

\newtheorem{Corollary2.3}[Theorem2.1]{\bf Corollary}
\begin{Corollary2.3}\label{Corollary2.3}
If $B$ is a $2\times 2$ matrix over a field $F$, then
\[\textnormal{Cen}_{M_2(F)}(B)=\left\{\begin{array}{l}M_2(F),\ \textit{ if }B\textit{ is a scalar matrix}\\\{f(B)\ |\ f(x)\in F[x]\},\textit{ if }B \textit{ is a non-scalar matrix}.\end{array}\right.\]
\end{Corollary2.3}
\noindent In this paper we denote the identity matrix by $E$.

\newtheorem{remCorollary2.3}[Theorem2.1]{\bf Remark}
\begin{remCorollary2.3}\label{remCorollary2.3}
\textnormal{Let $B=\left[\begin{array}{cc}e&f\\g&h\end{array}\right]\in M_2(R)$, $R$ a commutative ring. Elementary matrix multiplication shows that \begin{equation}\label{defA}A=\nolinebreak\left[\begin{array}{cc}a&b\\c&d\end{array}\right]\in\nolinebreak\textnormal{Cen}_{M_2(R)}(B)\end{equation} if and only if
\begin{equation}\label{basics}
(a-d)f=b(e-h),\quad bg=cf,\quad c(e-h)=(a-d)g
\end{equation}
if and only if $A'+vE$ and $B$ commute if and only if $A'+vE$ and $B'+wE$ commute if and only if $A'$ and $B'$ commute, where
\begin{equation}\label{defA'}
A'=\left[\begin{array}{cc}a-d&b\\c&0\end{array}\right]\qquad\textnormal{and}\qquad B'=\left[\begin{array}{cc}e-h&f\\g&0\end{array}\right].
\end{equation}}
\end{remCorollary2.3}
\vskip 0.5cm
\noindent Throughout the sequel, for $R$ a UFD and for a nonempty set $X\subset R$, we mean by $\textnormal{gcd}(X)$ an arbitrary greatest common divisor of $X$ in $R$. 

The following result is an extension of Corollary \ref{Corollary2.3} to UFD's.

\newtheorem{Corollary2.5b}[Theorem2.1]{\bf Corollary}
\begin{Corollary2.5b}\label{Corollary2.5b}
Let $B=\left[\begin{array}{cc}e&f\\g&h\end{array}\right]\in M_2(R)$, $R$ a UFD. Then
$\textnormal{Cen}_{M_2(R)}(B)$
\[=\left\{\begin{array}{l}(i)\:M_2(R),\ \textit{if $e=h$, $f=0$ and $g=0$ (i.e.~$B$ is a
scalar matrix)}\\\\
(ii)\left.\left\{m^{-1}w\left[\begin{array}{cc}e-h&f\\g&0\end{array}\right]+vE\right|v,w\in
 R\right\}, \begin{array}{l}\textit{if at least one}\\\textit{of $e-h,f,g$ is nonzero,}\end{array}\end{array}\right.\]
 where $m^{-1}$ is the inverse of $m:=\textnormal{gcd}(e-h,f,g)$ in the quotient field of~$R$.
\end{Corollary2.5b}
\begin{proof}
(ii) Suppose that at least one of $e-h$, $f$ and $g$ is nonzero. Let $A'$ and $B'$ be as in (\ref{defA'}). By the symmetry of the system of equations in (\ref{basics}) we may assume that $e-h\neq 0$. Then, using (\ref{basics}), $e-h|(a-d)f$ and $e-h|(a-d)g$ imply that~$e-h|m(a-d)$. Let $w\in R$ such that $m(a-d)=w(e-h)$. Then, again using~(\ref{basics}),~$(a-d)f=b(e-h)$ and $c(e-h)=(a-d)g$ imply that $mb=wf$ and~$mc=wg$. Thus $mA'=wB'$ and the result follows from Remark \ref{remCorollary2.3}.
\end{proof}

\newtheorem{Example2.5b}[Theorem2.1]{\bf Example}
\begin{Example2.5b}\label{Example2.5b}
\textnormal{Let $R$ be the UFD $\mathbb Z$ of integers, and let $B\nolinebreak=\nolinebreak\left[\begin{array}{cc}8&3\\6&2\end{array}\right].$ 
It follows from Corollary~\ref{Corollary2.5b}(ii) that
\begin{equation*}
\textnormal{Cen}_{M_2(\mathbb{Z})}(B)
=\left.\left\{\left[\begin{array}{cc}2w+v&w\\2w&v\end{array}\right]\right|v,w\in\mathbb{Z}\right\}.
\end{equation*}}
\end{Example2.5b}\vskip 0.5cm

For the remaining results in this section, let $\theta:R\to S$ be a ring
epimorphism and $\Theta:M_n(R)\to M_n(S)$ the induced epimorphism,
i.e.~$\Theta([b_{ij}])=[\theta(b_{ij})]$. We denote the annihilator of an element $r$ in a commutative ring $R$ by ann$_R(r)$. For the sake of notation, we will sometimes denote $\theta(b)$ by $\hat b$ and $\Theta(B)$ by $\widehat B$. Also, if there is no ambiguity, we simply write Cen$(B)$ instead of Cen$_{M_2(R)}(B)$ and Cen$(\widehat B)$ instead of Cen$_{M_2(S)}(\widehat B)$ for $B\in M_2(R)$, as well as ann$(\hat r)$ instead of ann$_{S}(\hat r)$ for~$r\in R$. If~$r\in R$ and $A\subseteq R$, then $rA$ denotes the set $\{ ra\ |\ a\in A\}$.

Throughout this paper and in particular in Section~4 we use the notation
\[\left[\begin{array}{cc}\mathcal B&\mathcal C\\\mathcal D&\mathcal E\end{array}\right]\quad\textrm{
to denote the set}\quad
\left\{\left.\left[\begin{array}{cc}b&c\\d&e\end{array}\right]\ \right|\ b\in\mathcal B, c\in\mathcal C, d\in \mathcal D, e\in \mathcal E\right\},\]
where $\mathcal B$, $\mathcal C$, $\mathcal D$ and $\mathcal E$ are subsets of a ring $R$.

The following result is straightforward.

\newtheorem{newlemma}[Theorem2.1]{\bf Lemma}
\begin{newlemma}\label{newlemma}
Let $S$ be a subring of a ring $T$ and let $s\in S$.
Then
\[\textnormal{Cen}_S(s)=S\cap \textnormal{Cen}_T(s).\]
\end{newlemma}

The following result is the inspiration behind Section 4.

\newtheorem{Lemma5}[Theorem2.1]{\bf Proposition}
\begin{Lemma5}\label{Lemma5}
Let $R$ be a commutative ring and let $B=[b_{ij}]\in M_n(R)$. Then
\[\Theta(\textnormal{Cen}(B))+[\mathcal{A}_{ij}]\subseteq
\textnormal{Cen}(\widehat B),\] where
\[\mathcal{A}_{ij}=\left(\displaystyle\bigcap_{k,\ k\neq j} \textnormal{ann}(\hat b_{jk})\right)\bigcap\left(\displaystyle\bigcap_{k,\ k\neq i} \textnormal{ann}(\hat b_{ki})\right)\bigcap\ \textnormal{ann}(\hat b_{ii}-\hat b_{jj}).\]
\end{Lemma5}
\begin{proof}
It follows easily that
\begin{equation}\label{eq2}
\Theta(\textnormal{Cen}(B))\subseteq \textnormal{Cen}(\widehat B).
\end{equation}
Now we show that
\begin{equation}\label{eq3}
[\mathcal{A}_{ij}]\subseteq \textnormal{Cen}(\widehat B).
\end{equation}
Let $[\hat
a_{ij}]\in[\mathcal{A}_{ij}]$. It follows that position $(r,t)$
of $\widehat B[\hat a_{ij}]-[\hat a_{ij}]\widehat B$ is equal to
\begin{equation*}
\hat b_{r1}\hat a_{1t}+\cdots +\hat b_{r,r-1}\hat a_{r-1,t}+\hat
b_{rr}\hat a_{rt}+\hat b_{r,r+1}\hat
a_{r+1,t}+\cdots+\hat b_{rn}\hat a_{nt}
-\end{equation*}\begin{equation}
(\hat a_{r1}\hat b_{1t}+\hat a_{r2}\hat b_{2t}+\cdots +\hat
a_{r,t-1}\hat b_{t-1,t}+\hat a_{rt}\hat b_{tt}+\hat a_{r,t+1}\hat
b_{t+1,t}+\cdots+\hat a_{rn}\hat b_{nt}).\qquad\label{eq4}
\end{equation}
Since $\hat a_{lt}\in \textnormal{ann}(\hat b_{rl})$ for every $l$
such that $l\neq r$, and $\hat a_{rq}\in \textnormal{ann}(\hat b_{qt})$ for every $q$ such that
$q\neq t$, according to the definition of
$[\mathcal{A}_{ij}]$, it follows that (\ref{eq4}) is equal to
\begin{equation}\label{eq5}
\hat b_{rr}\hat a_{rt}-\hat a_{rt}\hat b_{tt}=\hat a_{rt}(\hat
b_{rr}-\hat b_{tt}).
\end{equation}
Since $\hat a_{rt}\in \textnormal{ann}(\hat b_{rr}-\hat b_{tt})$, according to
the definition of $[\mathcal{A}_{ij}]$, it follows that (\ref{eq5})
is equal to $\hat 0$. Thus position $(r,t)$ of $[\hat
a_{ij}]\widehat B-\widehat B[\hat a_{ij}]$ is $\hat 0$. This proves
(\ref{eq3}).
\end{proof}

\section{$I$-invertibility in $R/I$ and $I$-matrices in $M_2(R/I)$, $R$ a UFD}
\noindent From here onwards, unless stated otherwise, we assume that $R$ is a UFD, $I$ is a nonzero ideal in $R$ and~$k:=\textnormal{gcd}(I)\neq 0$.
Let $\theta_I:R\to R/I$ and $\Theta_I:M_2(R)\to
M_2(R/I)$ be the
natural epimorphism and induced epimorphism respectively. We denote the image $\theta_I(b)$ of~$b\in R$ by $\hat b_I$ and the image $\Theta_I(B)$ of $B\in M_2(R)$ by $\widehat B_I$. However, if there is no ambiguity, then we simply write $\theta$, $\Theta$,~$\hat b$ and $\widehat B$ respectively.

The following results are trivial.

\newtheorem{Lemma8}[Theorem2.1]{\bf Lemma}
\begin{Lemma8}\label{Lemma8}
Let $R$ be a UFD. Then an element $\hat b=\theta(b)\in R/I$ is a zero
divisor if~\textnormal{gcd}$(b,k)\neq 1$.
\end{Lemma8}

\newtheorem{Lemma8a}[Theorem2.1]{\bf Lemma}
\begin{Lemma8a}\label{Lemma8a}
Let $R$ be a PID.
Then an element $\hat b\in R/\langle k\rangle$, $k\in R$, is invertible if and only if~\textnormal{gcd}$(b,k)=1$.
\end{Lemma8a}

\newtheorem{Def8}[Theorem2.1]{\bf Definition}
\begin{Def8}\label{Def8}
An $I$-pre-image of an element $\hat b\in R/I$ is a pre-image of~$\hat b$ in~$R$ of the form $r\delta$, where $\textnormal{gcd}(r,k)=1$ and ($\delta=0$ or $\delta|k$). If $\hat b=\hat 0$ we define~$\delta:=0$. We call $r$ and $\delta$ the relative prime part and divisor part of $r\delta$ respectively. We call~$\hat b$~$I$-invertible if $\hat r$ is invertible in~$R/I$ for at least one $I$-pre-image $r\delta$ of $\hat b$.
\end{Def8}

\newtheorem{Rem8}[Theorem2.1]{\bf Remark}
\begin{Rem8}\label{Rem8}\textnormal{It follows from Definition \ref{Def8} that if an element $\hat 0\neq\hat b\in R/I$ is~$I$-invertible, then there exists a $\hat c\in R/I$ such that $\hat c\hat b$ has a pre-image $\delta\in R$ which is a divisor of $k$.}\end{Rem8}%The following lemma is trivial to prove.

The converse of the above remark is not in general true. Here follows a counter example. 

\newtheorem{Example2c}[Theorem2.1]{\bf Example}
\begin{Example2c}\label{Example2c}
\textnormal{Let $R=\mathbb Z[x]$, let $I=\langle 5x^2\rangle$ and let $\hat b_I=\widehat{3x^2}$, then~$\hat b_{I}$ is not $I$-invertible, but \[\widehat{2}_I\hat b_I=\theta_I(6x^2-5x^2)=\widehat{x^2}_I.\]}
\end{Example2c}

We define the ideal $\delta^{-1}I:=\{\delta^{-1}a|a\in I\}\subset R$. The following result can be easily proved.

\newtheorem{Lemma8a'}[Theorem2.1]{\bf Lemma}
\begin{Lemma8a'}\label{Lemma8a'}
Let $\delta$ be the divisor part of an $I$-pre-image of $\hat 0\neq\hat b_I\in R/I$. There exists a $\hat c_I\in R/I$ such that $\hat c_I\hat b_I=\hat \delta_I$ if and only if $\widehat {b\delta^{-1}}_{\delta^{-1}I}$ is invertible in $R/\delta^{-1}I$, with inverse $\hat c_{\delta^{-1}I}$.
\end{Lemma8a'}

\newtheorem{Lemma8a''}[Theorem2.1]{\bf Lemma}
\begin{Lemma8a''}\label{Lemma8a''}
An element $\hat b_I\in R/I$ is $I$-invertible if and only if there exists an invertible element $\hat c\in R/I$ such that $\hat c_I\hat b_I=\hat\delta_I$, where $\delta$ is a divisor part of an~$I$-pre-image of $\hat b_I$.
\end{Lemma8a''}
\begin{proof}
If $\hat b_I$ is $I$-invertible then it follows directly from Definition \ref{Def8} that there exists an invertible element $\hat c_I$ such that $\hat c_I\hat b_I=\hat\delta_I$. Conversely, suppose there exists an invertible element $\hat c_I\in R/I$ such that $\hat c_I\hat b_I=\hat\delta_I$. Since $\hat c_I$ is invertible we have that $\hat b_I=\hat c_I^{-1}\hat\delta_I$. Let $c'\in R$ be a pre-image of $\hat c_I^{-1}$. Since $\hat c_I^{-1}$ is not a zero divisor it follows from Lemma \ref{Lemma8} that $\textnormal{gcd}(c',k)=1$. Since $c'\delta$ is an $I$-pre-image of $\hat b_I$ we have the desired result.
\end{proof}

The proof of the next result is constructive.

\newtheorem{Lemma10}[Theorem2.1]{\bf Lemma}
\begin{Lemma10}\label{Lemma10}
Every element in $R/I$ has an $I$-pre-image.
\end{Lemma10}

\begin{proof}
Let $\hat b\in R/I$. If $k$ is a unit, then the result follows trivially. Thus suppose $k$ is a nonzero nonunit. Since $R$ is a UFD there exist different primes~$p_1,\ldots,p_s$ such that
$k=p_1^{m_1}\cdots p_s^{m_s}$, where $m_1,\ldots,m_s\ge 1$. Since $1\cdot 0$ is an $I$-pre-image of $\hat 0$, suppose $\hat b$ is nonzero. Let $b$ be a pre-image of $\hat b$ in $R$. Again, because $R$ is a
UFD, $b$ can be expressed as $r_0p_1^{q_1}\cdots p_s^{q_s}$, where
$p_i\nmid\nolinebreak r_0$, for $i=1,\ldots,s$, and $q_1,\ldots,q_s\ge0.$
Therefore \textnormal{gcd}$(r_0,k)=1$, and
\[\hat b=\hat r_0\widehat{p_1^{q_1}}\cdots\widehat{p_s^{q_s}}.\]

Suppose we can show that each $\widehat{p_i^{q_i}}$ has a pre-image
$r_i\cdot p_i^{t_i}$, where \textnormal{gcd}$(r_i,k)=1$ and $t_i\le m_i$. Then we have that
\[
\hat b=\hat r_0\widehat
{(r_1p_1^{t_1})}\cdots\widehat{(r_sp_s^{t_{s}})}
=\hat r_0\hat r_1\cdots\hat
r_s\widehat{(p_1^{t_1}\cdots p_s^{t_s})}
=\theta({r{p_1^{t_1}\cdots p_s^{t_s}}}),
\]
where $r=r_0r_1\cdots r_s$. Since \textnormal{gcd}$(r_i,k)=1$ for $i=0,1,\ldots,s$,
 it follows that \textnormal{gcd}$(r,k)=\nolinebreak 1.$ Also, since $t_i\le m_i$ for $i=1,2\ldots,s,$ we have that
\[\delta:=p_1^{t_1}\cdots p_s^{t_s}|\underbrace{p_1^{m_1}\cdots
p_s^{m_s}}_{=k},\]
implying that $r\cdot \delta$ is an $I$-pre-image of $\hat b$ with relative prime part $r$ and divisor part~$\delta$.

Let us now prove that each $\widehat{p_i^{q_i}}$ has a pre-image
$r_i\cdot p_i^{t_i}$, where \textnormal{gcd}$(r_i,k)=1$ and $t_i\le m_i$.

If $q_i\le m_i$ then $p_i^{q_i}=1\cdot p_i^{q_i}$, where $t_i=q_i\le m_i$
and \textnormal{gcd}$(r_i,k)=1$, with $r_i=1$.
Thus we have the desired result.

Next we consider the case when $m_i<q_i$. Because $p_i^{m_i+1}\nmid k$, it follows that there exist an $a=a'k\in I$ such that $p_i\nmid a'$. Now since
\begin{eqnarray*}
\widehat{p_i^{q_i}}=\widehat{p_i^{q_i}+a'k}
\end{eqnarray*}
and
\begin{eqnarray*}
p_i^{q_i}+a'k=p_i^{q_i}+a'p_1^{m_1}\cdots p_s^{m_s}=p_i^{m_i}(p_i^{q_i-m_i}+a'p_1^{m_1}\cdots
p_{i-1}^{m_{i-1}}p_{i+1}^{m_{i+1}}\cdots p_s^{m_s}),
\end{eqnarray*}
it follows that $p_i^{m_i}\cdot r_i=r_i\cdot p_i^{m_i}$ is a pre-image of
$\widehat{p_i^{q_i}}$, where
\[r_i=p_i^{q_i-m_i}+a'p_1^{m_1}\cdots
p_{i-1}^{m_{i-1}}p_{i+1}^{m_{i+1}}\cdots p_s^{m_s}.\] Since
\[p_i|p_i^{q_i-m_i}(q_i>m_i)\qquad\textnormal{and}\qquad p_i\nmid a'p_1^{m_1}\cdots
p_{i-1}^{m_{i-1}}p_{i+1}^{m_{i+1}}\cdots p_s^{m_s},\] we have that
$p_i\nmid r_i$. Furthermore, for all
$l\in\{1,\ldots,i-1,i+1,\ldots,s\}$ it follows that \[p_l\nmid
p_i^{q_i-m_i}\qquad\textnormal{and}\qquad p_l|a'p_1^{m_1}\cdots
p_{i-1}^{m_{i-1}}p_{i+1}^{m_{i+1}}\cdots p_s^{m_s}\] implying that
$p_l\nmid r_i$. Thus $r_i$ and $k$ are relatively prime and
$t_i=m_i\le m_i$.
\end{proof}

We will now focus on the $I$-invertibility of elements in $R/I$.

The next result follows directly from Lemma \ref{Lemma8a}, Definition \ref{Def8} and Lemma~\ref{Lemma10}.

\newtheorem{Corollary1}[Theorem2.1]{\bf Corollary}
\begin{Corollary1}\label{Corollary1}
If $R$ is a PID, then every element in $R/I$ is $I$-invertible.
\end{Corollary1}

The next example illustrates the constructive proof of Lemma~\ref{Lemma10}.

\newtheorem{Example2}[Theorem2.1]{\bf Example}
\begin{Example2}\label{Example2}
\textnormal{Let $R=\mathbb{Z}$ and let $I=\langle 12\rangle$. Since $12=2^2\cdot3$  using the procedure in the proof of Lemma \ref{Lemma10}, it follows that
\begin{itemize}
\item[(a)]$\hat9_I=\theta_I({2^0\cdot3^2})=\theta_I({1\cdot(3^2+12)})=\theta_I({3(7)})=(\widehat{7\cdot3})_I$, where
\textnormal{gcd}$(7,12)\linebreak=\nolinebreak1$ and $3|12$.
%\item[(2)]$\hat 6_{10}=\theta_{10}(3\cdot 2\cdot 5^0)=(\widehat{3\cdot 2})_{10}$, where \textnormal{gcd}$(3,10)=1$ and $2|10$.
Since $\hat{7}_I$
is invertible in $\mathbb{Z}_{12}$ it follows that $\hat 9_I$ is $I$-invertible, as expected from Corollary \ref{Corollary1}.\newline
\end{itemize}
\indent Now, let $R=\mathbb Z[x]$ and let $I$ be a nonzero, not necessarily finite, ideal, with $2^4x^4\in I$ and $k:=\textnormal{gcd}(I)=2^3x^3$.
\begin{itemize}
\item[(b)] $\widehat{24x^5}_I+\widehat{8x^4}_I+\widehat{4x^2}_I=\theta_I(24x^5+8x^4+4x^2)=\theta_I((6x^3+2x^2+1)2^2x^2)$, where gcd$(6x^3+2x^2+1,2^3x^3)=1$ and $2^2x^2|2^3x^3$. Since $\widehat{6x^3}_I+\widehat{2x^2}_I+\hat1_I=\theta_I((3x+1)2x^2+1)$ is invertible in $R/I$ by Lemma \ref{Defnearlyprincipalrem}, $\widehat{24x^5}_I+\widehat{8x^4}_I+\widehat{4x^2}_I$ is $I$-invertible.
\end{itemize}}
\end{Example2}\vskip 0.15 cm

We already know from Example \ref{Example2c} that %will show in Example \ref{Example2b} that 
Corollary \ref{Corollary1} does not hold for $R$ a UFD in general, not even for the case when $I$ is a principal ideal. 

Lemma \ref{Corollary25}, Proposition \ref{Lemma23}, Remark \ref{Remark2b} and Lemma \ref{Defnearlyprincipalrem} will help us to determine when an element in~$R/I$ is not $I$-invertible in case $R$ is a UFD which is not a PID. In order to conclude that an element $\hat b\in R/I$ is not $I$-invertible (using Definition~\ref{Def8}), we have to show, for every $I$-pre-image $r\delta$ of $\hat b$, that $\hat r$ is not invertible in $R/I$. However, if $\hat b$ is principal (Definition~\ref{Defnearlyprincipal}), then we will show in Proposition \ref{Lemma23} that it suffices to show that $\hat r$ is not invertible in $R/I$ for at least one $I$-pre-image~$r\delta$ of~$\hat b$.

We first give a characterization of and establish a relationship between the divisor parts of the $I$-pre-images of an element in $R/I$.

\newtheorem{Lemma22}[Theorem2.1]{\bf Lemma}
\begin{Lemma22}\label{Lemma22}
Let $R$ be a UFD and let $\hat 0\neq\hat b\in R/I$. Then $\delta$ is a divisor part of an~$I$-pre-image of $\hat b$ if and only if $\textnormal{gcd}(b,k)=\delta$, i.e.~the divisor parts of the~$I$-pre-images of $\hat b$ are associates.
\end{Lemma22}
\begin{proof}
Let $r\delta$ be an $I$-pre-image of $\hat b$. Then $b=r\delta+sk$ for some~$s\in R$. Now, since $\textnormal{gcd}(r,k)=1$, it follows that $\textnormal{gcd}(b,k)=\textnormal{gcd}(r\delta+sk,k)=\textnormal{gcd}(\delta,k)=\delta$.

For the converse, note that since all the greatest common divisors of $b$ and $k$ are associates and every element in $R/I$ has at least one $I$-pre-image, by Lemma~\ref{Lemma10}, the result will follow if we can show that for an arbitrary unit $t$, $t\delta$ is also a divisor part of some $I$-pre-image of $\hat b$. Since $\widehat{rt^{-1}t\delta}=\widehat{r\delta}=\hat b$, $\textnormal{gcd}(rt^{-1},k)=1$ and $t\delta|k$, the result follows.
\end{proof}

The following result follows trivially from Lemma \ref{Lemma22}.

\newtheorem{Corollary25}[Theorem2.1]{\bf Lemma}
\begin{Corollary25}\label{Corollary25}
Let $\hat0\neq\hat b\in R/I$.
If $\textnormal{gcd}(b,k)=1$, then $\hat b$ is $I$-invertible if and only if $\hat b$ is invertible in $R/I$.
\end{Corollary25}

\newtheorem{Remark25}[Theorem2.1]{\bf Remark}
\begin{Remark25}\label{Remark25}
\textnormal{Note that if $k$ is a unit, it follows from Lemma \ref{Corollary25} that every~$\hat0\neq\hat b\in R/I$ is $I$-invertible if and only if $\hat b$ is invertible in $R/I$.}
\end{Remark25}

\newtheorem{Defnearlyprincipal}[Theorem2.1]{\bf Definition}
\begin{Defnearlyprincipal}\label{Defnearlyprincipal}
Let $R$ be a UFD, let $k=p_1^{m_1}\cdots p_s^{m_s}\in R$ be a nonunit, with $p_1,\ldots, p_s$ different primes and $m_1,\ldots,m_s\ge1$, and let $\hat b\in R/I$. If~$\delta:=\textrm{gcd}(b,k)=p_1^{q_1}\cdots p_s^{q_s}$, where $0\le q_i<m_i$ for $i=1,\ldots,s$, then we call $\hat b$ a principal element of $R/I$. 
If $\widehat{\delta^{-1}k}$ is principal, i.e.~$\delta=p_1^{q_1}\cdots p_s^{q_s}$, where $q_i\ge1$ for~$i=1,\ldots,s$, we call $\hat b$ $q$-principal.
\end{Defnearlyprincipal}

\newtheorem{Lemma23}[Theorem2.1]{\bf Proposition}
\begin{Lemma23}\label{Lemma23}
Let $R$ be a UFD, $k$ be a nonunit and let $\hat 0\neq\hat b\in\nolinebreak R/I$ be principal, then either $\hat r$ is invertible in $R/I$ for every $I$-pre-image $r\delta$ of~$\hat b$ or no such~$\hat r$ is invertible in~$R/I$.
\end{Lemma23}

\begin{proof}
Since, according to Lemma~\ref{Lemma10}, there exists a pre-image $r\delta$ of $\hat b$ in $R$, with \textnormal{gcd}$(r,k)=1$, all the pre-images, and in particular all the $I$-pre-images, of $r\delta$ are of the form
\begin{equation}\label{eq23.1}
r\delta+cp_1^{m_1}p_2^{m_2}\cdots p_s^{m_s},
\end{equation}
where $cp_1^{m_1}p_2^{m_2}\cdots p_s^{m_s}\in I$. Because, according to Lemma \ref{Lemma22}, the divisor parts of all the~$I$-pre-images of $\hat b$ are of the form $u\delta$, where $u$ is a unit in $R$, it follows from~(\ref{eq23.1}) that the relative prime parts of all the $I$-pre-images of $\hat b$ are of the form
\begin{equation}
u^{-1}r+cu^{-1}p_1^{m_1-q_1}\cdots p_s^{n_s-q_s},
\end{equation}
where $cp_1^{m_1}p_2^{m_2}\cdots p_s^{m_s}\in I$ and $u\in R$ is a unit.

Now, suppose $\hat r$ is invertible in $R/I$ with inverse $\hat y$. In other words
\[yr=1+dp_1^{m_1}p_2^{m_2}\cdots p_s^{m_s},\]where $dp_1^{m_1}p_2^{m_2}\cdots p_s^{m_s}\in I$. If we can show that the image under $\theta$ of the relative prime part of an arbitrary $I$-pre-image of $\hat b$ is invertible, then we are finished.

Let $u^{-1}r+cu^{-1}p_1^{m_1-q_1}p_2^{m_2-q_2}\cdots p_s^{m_s-q_s}$ be the relative prime part of an arbitrary $I$-pre-image of $\hat b$. Furthermore, let $l\in\mathbb{N}$ such that
\begin{equation}\label{>}
2^l>\max\left\{\left.\frac{m_i}{m_i-q_i}\ \right|\ i\in\{1,\ldots,s\}\right\}>0.
\end{equation}
For the sake of notation, let
\[v=dp_1^{q_1}\cdots p_s^{q_s}+cy\textrm{ and }w=p_1^{m_1-q_1}p_2^{m_2-q_2}\cdots p_s^{m_s-q_s}.\] Then
\begin{eqnarray*}
&&(u^{-1}r+cu^{-1}p_1^{m_1-q_1}p_2^{m_2-q_2}\cdots p_s^{m_s-q_s})yu(1-vw)(1+(vw)^{2^1})\\
&&(1+(vw)^{2^{2}})(1+(vw)^{2^{3}})(1+(vw)^{2^{4}})\cdots(1+(vw)^{2^{l-1}})\\
&=&(1+dp_1^{m_1}p_2^{m_2}\cdots p_s^{m_s}+cyp_1^{m_1-q_1}p_2^{m_2-q_2}\cdots p_s^{m_s-q_s})(1-vw)\\&&(1+(vw)^{2^{1}})
(1+(vw)^{2^{2}})\cdots(1+(vw)^{2^{l-1}})\\
&=&(1+vw)(1-vw)(1+(vw)^{2^{1}})\cdots(1+(vw)^{2^{{l-1}}})\\
&=&1-(vw)^{2^{l}}.
\end{eqnarray*}
Let $1\le i\le s$. Since $m_i>q_i$, it follows from (\ref{>}) that
\[2^{l}(m_i-q_i)>\frac{m_i}{m_i-q_i}(m_i-q_i)=m_i,\]
and so
\[w^{2^{l}}
=ap_1^{m_1}p_2^{m_2}\cdots p_s^{m_s}\]
for some $a\in R$. Since $dp_1^{m_1}\cdots p_s^{m_s}\in I$ and $cp_1^{m_1}\cdots p_s^{m_s}\in I$ imply that $vw^{2^l}\in I$, it follows that $(vw)^{2^{l}}\in I$. Therefore
\begin{eqnarray*}
&&\theta\left((u^{-1}r+cu^{-1}p_1^{m_1-q_1}p_2^{m_2-q_2}\cdots p_s^{m_s-q_s})yu(1-vw)\right.\\
&&\left.(1+(vw)^{2^{1}})(1+(vw)^{2^{2}})\cdots(1+(vw)^{2^{l-1}})\right)\\
&=&\theta\left(1-(vw)^{2^l}\right)\\
&=&\hat 1.
\end{eqnarray*}
Hence, we conclude that
\[\theta\left(yu(1-vw)(1+(vw)^{2^{1}})(1+(vw)^{2^{2}})\cdots(1+(vw)^{2^{l-1}})\right)\]
is the inverse of the image under $\theta$ of the relative prime part of the arbitrary chosen~$I$-pre-image of $\hat b$.
\end{proof}

\newtheorem{Remark2b}[Theorem2.1]{\bf Remark}
\begin{Remark2b}\label{Remark2b}
\textnormal{
Note that if $I=\langle p^n\rangle$, for a prime $p\in R$ and $n>0$, then every $\hat 0\neq\hat b\in R/I$ is principal. Thus Proposition \ref{Lemma23} is applicable to all nonzero elements in $R/I$. Furthermore, it is helpful to notice that every pre-image of $\hat 0\neq\hat b$ is an~$I$-pre-image.}
\end{Remark2b}

Next we show that Proposition \ref{Lemma23} does not hold in general if $q_i=m_i$ for some~$i$.

\newtheorem{Example23}[Theorem2.1]{\bf Example}
\begin{Example23}\label{Example23}
\textnormal{Let $R=\mathbb Z[x]$, let $k=2x$ (with $2$ and $x$ primes in $\mathbb Z[x]$) and let $I=\langle 2x\rangle$. Consider $\hat 0\neq\hat x\in\mathbb Z[x]/\langle 2x\rangle$. Then $1\cdot x$ and $3\cdot x$ are $\langle 2x\rangle$-pre-images of $\hat x$ with relative prime parts $1$ and $3$ respectively, and $\hat 1$ is invertible in $\mathbb Z[x]/\langle 2x\rangle$, but $\hat 3$ is not.}
\end{Example23}

\newtheorem{Lemma2e}[Theorem2.1]{\bf Lemma}
\begin{Lemma2e}\label{Lemma2e}
Let $k$ be a nonunit and let $\hat 0_I\neq\hat b_I\in R/I$ be principal with $\delta$ the divisor part of an~$I$-pre-image of $\hat b$. Then there exists a $\hat c_I\in R/I$ such that~$\hat c_I\hat b_I=\hat\delta_I$ if and only if $\hat b_I$ is $I$-invertible.
\end{Lemma2e}
\begin{proof}
Let $r\delta$ be an $I$-pre-image of $\hat b_I$ and suppose that there exists a $\hat c_I\in R/I$ such that $\hat c_I\hat b_I=\hat\delta_I$. Then it follows from Lemma \ref{Lemma8a'} that $cr=cb\delta^{-1}=1+\gamma\delta^{-1}k$, for some $\gamma\in R$ such that $\gamma k\in I$. Suppose $k=p_1^{q_1}\cdots p_s^{q_s}$, for $p_1,\ldots,p_s$ prime and $q_1,\ldots, q_s\ge 1$. Since $\hat b_I$ is principal, it follows that $k\delta^{-1}$ is of the form $w=p_1^{v_1}\ldots p_s^{v_s}$, where $v_1,\ldots,v_s\ge 1$. Now let $l\in \mathbb N$ such that $2^l>\max\{q_1,\ldots,q_s\}$. Then
\[(1+\gamma\delta^{-1}k)(1-\gamma\delta^{-1}k)(1+(\gamma\delta^{-1}k)^{2^1})\cdots(1+(\gamma\delta^{-1}k)^{2^{l-1}})=1-(\gamma\delta^{-1}k)^{2^{l}}\]which implies that
\[cr(1-\gamma\delta^{-1}k)(1+(\gamma\delta^{-1}k)^{2^1})\cdots(1+(\gamma\delta^{-1}k)^{2^{l-1}})=1-(\gamma\delta^{-1}k)^{2^{l}},\]
where $k|(\delta^{-1}k)^{2^l}$. Since $\gamma k\in I$, it follows that $(\gamma\delta^{-1}k)^{2^{l}}\in I$. Hence $\hat r_I$ is invertible in $R/I$ and we can conclude that $\hat b_I$ is $I$-invertible in $R/I$. The converse follows from Remark \ref{Rem8}.
\end{proof}

\newtheorem{principal}[Theorem2.1]{\bf Remark}
\begin{principal}\label{principal}
\textnormal{Let $k$ be a nonunit and let $\hat 0_I\neq\hat b_I\in R/I$ be principal. Using, Lemma \ref{Lemma8a''} and Lemma \ref{Lemma2e} it is only necessary to consider invertible elements in~$R/I$ to determine whether there exists a $\hat c_I$ in $R/I$ such that $\hat b_I\hat c_I=\hat\delta_I$, where $\delta$ is a divisor part of an $I$-pre-image of $\hat b_I$.}
\end{principal}

The following result will help us to determine whether an image of a relative prime part of an $I$-pre-image of an element is invertible in $R/I$ and can be proved by a similar method than the method in the proof of Lemma \ref{Lemma2e}.

\newtheorem{Defnearlyprincipalrem}[Theorem2.1]{\bf Lemma}
\begin{Defnearlyprincipalrem}\label{Defnearlyprincipalrem}
{Let $k\in R$ be a nonzero nonunit. If $\hat b\in R/I$ has a pre-image of the form $b'+1$, where $\hat{b'}$ is $q$-principal and $\hat{b'}\hat k\in I$, then $\hat b$ is invertible in $R/I$ (see Example~\ref{Example2}(b)).}
\end{Defnearlyprincipalrem}

\newtheorem{Defnearlyprincipalremrem}[Theorem2.1]{\bf Remark}
\begin{Defnearlyprincipalremrem}\label{Defnearlyprincipalremrem}
\textnormal{The converse of Lemma \ref{Defnearlyprincipalrem} is not in general true. For example~$\hat 3=\widehat{(2+1)}$ is invertible in $\mathbb Z_5$, although $5\nmid 2$.}
\end{Defnearlyprincipalremrem}

\newtheorem{Example2b}[Theorem2.1]{\bf Example}
\begin{Example2b}\label{Example2b}
\textnormal{Let $R=F[x,y]$ and let $I:=\langle y^5\rangle$. Since
$\widehat{x^5}_I$ is not invertible in $F[x,y]/\langle
y^5\rangle$, we conclude from Remark \ref{Remark2b} that $\widehat{x^5}_I$ is not $I$-invertible. Because \textnormal{gcd}$(x^5,y^5)=1$ we could also concluded from Lemma~\ref{Corollary25} that $\widehat{x^5}_I$ is not~$I$-invertible.}
\end{Example2b}

\newtheorem{Def10a}[Theorem2.1]{\bf Definition}
\begin{Def10a}\label{Def10a}
We call a matrix $\left[\begin{array}{cc}\hat e_I&\hat f_I\\\hat g_I&\hat h_I\end{array}\right]\in M_2(R/I)$ an $I$-matrix \linebreak if $\langle \hat e_I-\hat h_I,\hat f_I\rangle=\langle\hat t_I\rangle$ or $\langle \hat e_I-\hat h_I, \hat g_I\rangle=\langle\hat t_I\rangle$ or $\langle \hat f_I, \hat g_I\rangle=\langle\hat t_I\rangle$, where $t|k$.
\end{Def10a}

The following result is easy to prove.

\newtheorem{Lemma1}[Theorem2.1]{\bf Lemma}
\begin{Lemma1}\label{Lemma1}
Let $\hat a_I,\hat b_I\in R/I$. If $\langle\hat a_I,\hat b_I\rangle=\langle \hat t_I\rangle$, where $t|k$, then $t=\textnormal{gcd}(a,b,k)$.
\end{Lemma1}

The following results can be used to determine whether a matrix is an $I$-matrix.

\newtheorem{kinv}[Theorem2.1]{\bf Lemma}
\begin{kinv}\label{kinv}
A matrix is an $I$-matrix if it satisfies the following conditions:
\begin{itemize}
\item[(i)]For at least one of the three elements $\hat e_I-\hat h_I,\ \hat f_I$ and $\hat g_I$, say $\hat\alpha_I$, there exists a $\hat c_I\in R/I$ such that $\hat c_I\hat\alpha_I=\hat\delta_I$, where $r\delta$ is an $I$-pre-image of $\hat\alpha_I$ that has divisor part $\delta$; pick such an element, and call the remaining two elements $\hat a_I$ and $\hat b_I$, say.
\item[(ii)]For at least one of the elements $\hat a_{\langle\delta\rangle}$ and $\hat b_{\langle\delta\rangle}$, say $\hat\beta_{\langle\delta\rangle}$, there exists a $\hat d_{\langle\delta\rangle}\in R/\langle\delta\rangle$ such that $\hat d_{\langle\delta\rangle}\hat\beta_{\langle\delta\rangle}=\hat t_{\langle\delta\rangle}$, where $t|\delta$.
\end{itemize}
\end{kinv}

\newtheorem{remkinv}[Theorem2.1]{\bf Remark}
\begin{remkinv}\label{remkinv}
\textnormal{Note that if Lemma \ref{kinv}(i) is satisfied, with $\delta$ a unit, then Lemma \ref{kinv}(ii) is always satisfied.}
\end{remkinv}

The following result is in some cases helpful to determine when a matrix is not an $I$-matrix.

\newtheorem{Lemma2}[Theorem2.1]{\bf Lemma}
\begin{Lemma2}\label{Lemma2}
Let $\hat a_I,\hat b_I\in R/I$ and suppose that there exists a $\hat c_I$ such that $\hat c_I\hat a_I=\hat\delta_I$, with $\delta$ a divisor part of an $I$-pre-image of $\hat a_I$. Then $\langle\hat a_I,\hat b_I\rangle=\langle \hat t_I\rangle$, where $t|k$, if and only if there exists a $\hat d_{\langle\delta\rangle}$ such that $\hat d_{\langle\delta\rangle}\hat b_{\langle\delta\rangle}=\hat t_{\langle\delta\rangle}$.
\end{Lemma2}
\begin{proof}
Suppose there exists a $\hat c_I\in R/I$ such that $\hat c_I\hat a_I=\hat\delta_I$, with $\delta$ a divisor part of an $I$-pre-image of $\hat a_I$. 

Using Lemma \ref{Lemma1}, suppose that $\langle\hat a_I,\hat b_I\rangle=\langle\hat\delta_I,\hat b_I\rangle=\langle\hat t_I\rangle$, where $t=\textnormal{gcd}(\delta,b,k)=\textnormal{gcd}(\delta,b)$. Then, since $t|\delta|k$, $\alpha \delta+\beta b\equiv t+I$, for some $\alpha,\beta\in R$, implies that $\alpha \delta+\beta b= t+\gamma\delta$, for some $\gamma\in R$, and so~$\beta b=t+(\gamma-\alpha)\delta$. The converse follows trivially. 
\end{proof}

\newtheorem{Lemma2a}[Theorem2.1]{\bf Lemma}
\begin{Lemma2a}\label{Lemma2a}
If $\hat e-\hat h,\ \hat f$ or $\hat g$ is invertible in $R/I$ then $\left[\begin{array}{cc}\hat e&\hat f\\\hat g&\hat h\end{array}\right]\in M_2(R/I)$ is an $I$-matrix.
\end{Lemma2a}

\begin{proof}
Suppose $\hat c_I\in\{\hat e_I-\hat h_I,\hat f_I,\hat g_I\}$ is invertible in $R/I$. Then it follows from Lemma \ref{Corollary25} that $\hat c_I$ is $I$-invertible with an $I$-pre-image $c\cdot 1$ that has divisor part~$1$, and so the result follows from Remark \ref{Rem8}, Lemma \ref{kinv} and Remark \ref{remkinv}.
\end{proof}

The following result follows directly from Corollary \ref{Corollary1}, Remark \ref{Rem8} and \linebreak Lemma~\ref{kinv}.

\newtheorem{Corollary2}[Theorem2.1]{\bf Corollary}
\begin{Corollary2}\label{Corollary2}
If $R$ is a PID, then every matrix in $M_2(R/I)$ is an $I$-matrix.
\end{Corollary2}

We show that Corollary \ref{Corollary2} does not hold for UFD's in general.

\newtheorem{Example3}[Theorem2.1]{\bf Example}
\begin{Example3}\label{Example}
\textnormal{Let $R=\mathbb Z[x]$ and let $I$ be a nonzero ideal in $R$, with $2^4x^4\in I$ and $k=2^3x^3$. We exhibit~(a)~ a matrix which is an $I$-matrix and (b) a matrix which is not an $I$-matrix.
\newline(a) Let
\begin{equation*}
\widehat B_I=\left[\begin{array}{cc}\widehat{7x^2}_I&\widehat{24x^5}_I+\widehat{8x^4}_I+\widehat{4x^2}_I\\\widehat{14x}_I&\hat 0_I\end{array}\right]\in M_2(R/I).
\end{equation*}
We have already seen in Example~\ref{Example2}(b) that $\widehat{24x^5}_I+\widehat{8x^4}_I+\widehat{4x^2}_I$ is $I$-invertible with divisor part $\delta=2^2x^2$. Since $\widehat{7x^2}_{\langle\delta\rangle}=\widehat{-1x^2}_{\langle\delta\rangle}$, it follows that $\widehat{7x^2}_{\langle\delta\rangle}$ is $\langle\delta\rangle$-invertible and therefore, using Remark \ref{Rem8} and Lemma~\ref{kinv}, $\widehat B_I$ is an $I$-matrix.
\newline (b) Let
\begin{equation*}
\widehat B_I=\left[\begin{array}{cc}\hat{3}_I&\widehat{24x^5}_I+\widehat{8x^4}_I+\widehat{4x^2}_I\\\widehat{14x}_I&\hat 0_I\end{array}\right]\in M_2(R/I).
\end{equation*}
We first consider the ideals $\langle\hat 3_I,\widehat{24x^5}_I+\widehat{8x^4}_I+\widehat{4x^2}_I\rangle$ and $\langle\widehat{14x}_I,\widehat{24x^5}_I+\widehat{8x^4}_I+\widehat{4x^2}_I\rangle$. We have already seen in Example~\ref{Example2}(b) that $\widehat{24x^5}_I+\widehat{8x^4}_I+\widehat{4x^2}_I$ is $I$-invertible with divisor part $\delta=2^2x^2$. Since~$\hat 3_{\langle\delta\rangle}$ and $\widehat{14x}_{\langle\delta\rangle}=\hat7_{\langle\delta\rangle}\widehat{2x}_{\langle\delta\rangle}$ are both principal, it follows from Proposition \ref{Lemma23} that~$\hat 3_{\langle\delta\rangle}$ and $\widehat{14x}_{\langle\delta\rangle}$ are both not~$\langle\delta\rangle$-invertible. Therefore it follows from Lemma~\ref{Lemma2e}, Lemma~\ref{Lemma2} and Lemma~\ref{Lemma1} that $\widehat B_I$ is an~$I$-matrix if and only if $\langle\hat 3_I,\widehat{14x}_I\rangle=R/I$. Since this is not the case~$\widehat B_I$ is a non-$I$-matrix.} 
\end{Example3}\vskip 0.3cm

\section{The centralizer of an $I$-matrix}

The purpose of this section is to obtain a concrete description of the centralizer of an $I$-matrix in $M_2(R/I)$, $R$ a UFD and a nonzero ideal $I$ in $R$, with $k:=\textnormal{gcd}(I)$, by showing that the converse containments $\supseteq$ hold in Proposition \ref{Lemma5}. We also provide an example of a UFD, which is not a PID, and a non-$I$-matrix in~$M_2(R/I)$ for which the mentioned converse containment does not hold. We conclude with an example where we show that if $R$ is a UFD and $R/I$ is such that $R/I$ is not an integral domain, then for every $n\ge 3$ there is a matrix in $M_n(R)$ for which we do not have equality in Proposition~\ref{Lemma5}. Note that we still assume that $\theta_I:R\to R/I$ and~$\Theta_I: M_2(R)\to M_2(R/I)$ are the natural and induced epimorphism respectively.

\newtheorem{Theorem21}[Theorem2.1]{\bf Theorem}
\begin{Theorem21}\label{Theorem21}
Let $R$ be a UFD, $I$ a nonzero ideal in $R$, and let $\widehat B_I=\left[\begin{array}{cc}\hat e_I&\hat f_I\\\hat g_I&\hat h_I\end{array}\right]\in M_2(R/I)$ be an $I$-matrix, then
\begin{eqnarray}
\qquad\textnormal{Cen}(\widehat B)&=&\Theta(\textnormal{Cen}(B))+\left[\begin{array}{cc}\hat 0&\textnormal{ann}(\hat g)\cap
\textnormal{ann}(\hat e-\hat h)\\\textnormal{ann}(\hat f)\cap
\textnormal{ann}(\hat e-\hat h)&\textnormal{ann}(\hat f)\cap \textnormal{ann}(\hat g)\end{array}\right]\label{eq70}\\
&=&\Theta(\textnormal{Cen}(B))+\left[\begin{array}{cc}\textnormal{ann}(\hat f)\cap \textnormal{ann}(\hat g)&\textnormal{ann}(\hat g)\cap
\textnormal{ann}(\hat e-\hat h)\\\textnormal{ann}(\hat f)\cap
\textnormal{ann}(\hat e-\hat h)&\hat 0\end{array}\right]\label{eq71}\\
&=&\Theta(\textnormal{Cen}(B))+\left[\begin{array}{cc}\textnormal{ann}(\hat f)\cap \textnormal{ann}(\hat g)&\textnormal{ann}(\hat g)\cap
\textnormal{ann}(\hat e-\hat h)\\\textnormal{ann}(\hat f)\cap
\textnormal{ann}(\hat e-\hat h)&\textnormal{ann}(\hat f)\cap \textnormal{ann}(\hat g)\end{array}\right].\label{eq72}
\end{eqnarray}
\end{Theorem21}
\begin{proof}
By the symmetry in (\ref{basics}) it is sufficient to consider the case where $\langle \hat f,\hat g\rangle=\langle \hat t\rangle$, $t=\textnormal{gcd}(f,g,k)$ by Lemma \ref{Lemma1}. Suppose $\hat A=\left[\begin{array}{cc}\hat a&\hat b\\\hat c&\hat d\end{array}\right]\in M_2(R/I)$ such that $\widehat A\widehat B=\widehat B\widehat A$, i.e.~$A\in M_2(R)$ such that $AB\equiv BA+I$. Since $c(e-h)\equiv (a-d)g+I$ and $t|g,k$ it follows that $t|c(e-h)$. Let $m=\textnormal{gcd}(e-h,f,g,k)$, then $\textnormal{gcd}(e-h,t)=m$, which implies that $t|cm$. Similarly $(a-d)f\equiv b(e-h)+I$ yields $t|bm$. Since $\langle \hat f,\hat g\rangle=\langle \hat t\rangle$, there exists an $\hat\alpha,\hat\beta\in R/I$ such that $\hat t=\hat\alpha \hat f+\hat\beta \hat g$, i.e.~$t\equiv\alpha f+\beta g +I$. Let $w\in R$ such that $w=\alpha b+\beta c$, then $t|wm$. Let $v\in R$ such that $vt=wm$. It follows from (\ref{basics}), using the notation of Remark \ref{remCorollary2.3}, that
\[fA'\equiv bB'+I\qquad gA'\equiv cB'+I\quad\textnormal{and so}\quad wB'\equiv tA'+I.\]
Write $B'$ as $mB''$, then $vtB''= wmB''=wB'\equiv tA'+I$. Let $\widehat K=\left[\begin{array}{cc}\hat{e'}&\hat{f'}\\\hat{g'}&\hat{h'}\end{array}\right]$ be the image of $A'-vB''$ in $M_2(R/I)$ and $L=vB''$, then $L\in\textnormal{Cen}(B)$, by Lemma \ref{Corollary2.5b},
\[\hat t\widehat K=\hat 0\qquad\textnormal{and}\qquad\widehat{A'}=\widehat L+\widehat K.\] Here $\widehat K$ commutes with $\widehat {B'}$, and hence with $\widehat B$, and therefore $(\hat{e'}-\hat{h'})\hat f=\hat{f'}(\hat e-\hat h)$, $\hat{f'}g=\hat{g'}\hat f$ and $\hat{g'}(\hat e-\hat h)=(\hat{e'}-\hat{h'})\hat g$. But $(\hat{e'}-\hat{h'})\hat f=\hat0$, since $(\hat{e'}-\hat{h'})\hat t=\hat 0$ and $t|f$. Similarly $(\hat{e'}-\hat{h'})\hat g=\hat0$, $\hat{f'}\hat g=\hat 0$, $\hat{f'}(\hat e-\hat h)=\hat0$, $\hat{g'}\hat f=\hat0$ and $\hat{g'}(\hat e-\hat h)=\hat0$. Hence
\[\widehat K\subseteq\left[\begin{array}{cc}\textnormal{ann}(\hat f)\cap \textnormal{ann}(\hat g)&\textnormal{ann}(\hat g)\cap
\textnormal{ann}(\hat e-\hat h)\\\textnormal{ann}(\hat f)\cap
\textnormal{ann}(\hat e-\hat h)&\textnormal{ann}(\hat f)\cap \textnormal{ann}(\hat g)\end{array}\right].\]
Since $\widehat{A}=\widehat{A'}+\hat d\widehat{E}$, we have the containment $\subseteq$ in (\ref{eq72}). The converse follows from Proposition \ref{Lemma5}.
\end{proof}

\newtheorem{Example19}[Theorem2.1]{\bf Example}
\begin{Example19}\label{Example19}
\textnormal{Consider $\widehat B_I=\left[\begin{array}{cc}\widehat{7x^2}_I&\widehat{24x^5}_I+\widehat{8x^4}_I+\widehat{4x^2}_I\\\widehat{14x}_I&\hat 0_I\end{array}\right]\in M_2(R/I)$ in Example~\ref{Example}(a), with $I=\langle 5\cdot 2^3x^3,2^4x^4\rangle$. We use Theorem~\ref{Theorem21}, (\ref{eq71}), to obtain $\textnormal{Cen}(\widehat B)$. According to Corollary \ref{Corollary2.5b}(ii)
\begin{equation}
\textnormal{Cen}(B)=\left\{\left[\left.\begin{array}{cc}h_1+7xh_2&(24x^4+8x^3+4x)h_2\\14h_2&h_1\end{array}\right]\right|h_1,h_2\in \mathbb Z[x]\right\}.\label{eq82}
\end{equation}
Furthermore, \begin{eqnarray*}\textnormal{ann}(\widehat{7x})\cap\textnormal{ann}(\widehat{14x})&=&\langle \widehat{5\cdot2^3x^2},\widehat{2^4x^3}\rangle,\\ 
\textnormal{ann}(\widehat{14x})\cap \textnormal{ann}(\widehat{24x^5}+\widehat{8x^4}+\widehat{4x^2})&=&\langle\widehat{5\cdot2^2x^2},\widehat{2^3x^3}\rangle\\\textnormal{and }\textnormal{ann}(\widehat{7x})\cap \textnormal{ann}(\widehat{24x^5}+\widehat{8x^4}+\widehat{4x^2})&=&\langle\widehat{5\cdot2^3x^2},\widehat{2^4x^3}\rangle\end{eqnarray*} and so it follows from (\ref{eq82}) and Theorem~\ref{Theorem21}, (\ref{eq70}), that
\begin{eqnarray*}
\textnormal{Cen}(\widehat B)
&=&\Theta\left(\left\{\left[\left.\begin{array}{cc}h_1+7xh_2&(24x^4+8x^3+4x)h_2\\14h_2&h_1\end{array}\right]\right|h_1,h_2\in \mathbb Z[x]\right\}\right)\\&&+\left[\begin{array}{cc}\hat0 &\langle \widehat{5\cdot2^3x^2},\widehat{2^4x^3}\rangle\\\langle\widehat{5\cdot2^3x^2},\widehat{2^4x^3}\rangle&\langle\widehat{5\cdot2^2x^2},\widehat{2^3x^3}\rangle\end{array}\right].\\
%&=&\left\{\left.\left[\begin{array}{cc}\hat h_1+\hat h_3\hat x&\hat x\hat h_2\\\hat x\hat h_2&\hat h_1-\hat y\hat h_2\end{array}\right]\right|\ \hat h_1,\hat h_2,\hat h_3\in F[x,y]/\langle x^2\rangle\right\}.
\end{eqnarray*}
}
\end{Example19}
\vskip 0.5cm
\newtheorem{Remark19}[Theorem2.1]{\bf Remark}
\begin{Remark19}\label{Remark19}
\textnormal{Note that in the above example \[\Theta(\textnormal{Cen}(B))\not\subseteq\left[\begin{array}{cc}\textnormal{ann}(\hat f)\cap \textnormal{ann}(\hat g)&\textnormal{ann}(\hat g)\cap \textnormal{ann}(\hat e-\hat h)\\\textnormal{ann}(\hat f)\cap \textnormal{ann}(\hat e-\hat h)&\textnormal{ann}(\hat f)\cap \textnormal{ann}(\hat g)\end{array}\right]\]and that
\[\left[\begin{array}{cc}\textnormal{ann}(\hat f)\cap \textnormal{ann}(\hat g)&\textnormal{ann}(\hat g)\cap \textnormal{ann}(\hat e-\hat h)\\\textnormal{ann}(\hat f)\cap \textnormal{ann}(\hat e-\hat h)&\textnormal{ann}(\hat f)\cap \textnormal{ann}(\hat g)\end{array}\right]\not\subseteq\Theta(\textnormal{Cen}(B)).
\]
}
\end{Remark19}\vskip 0.5cm

According to Corollary~\ref{Corollary2}, Theorem~\ref{Theorem21} applies to all $2\times2$ matrices over factor rings $R/I$, where $R$ is a PID. In other words, we have equality in Proposition~\ref{Lemma5} for all $2\times 2$ matrices over factor rings of PID's. This is not the case for all $2\times 2$ matrices over factor rings $R/I$, where $R$ is a UFD, as the following example shows.

\newtheorem{Example20}[Theorem2.1]{\bf Example}
\begin{Example20}\label{Example20}
\textnormal{Consider $\widehat B=\left[\begin{array}{cc}\hat x+\hat y&\hat y\\\hat x&\hat x\end{array}\right]\in M_2(F[x,y]/\langle x^2\rangle)$. By Corollary~\ref{Corollary2.5b}(ii)
\begin{equation}
\textnormal{Cen}(B)=\left\{\left.\left[\begin{array}{cc}h_1&yh_2\\xh_2&h_1-yh_2\end{array}\right]\right|\ h_1,h_2\in F[x,y]\right\}.\label{eq90}
\end{equation}
The second term in the righthand side of (\ref{eq71}) is
\[\left[\begin{array}{cc}\textnormal{ann}(\hat y)\cap \textnormal{ann}(\hat x)&\textnormal{ann}(\hat x)\cap \textnormal{ann}(\hat y)\\ \textnormal{ann}(\hat y)\cap \textnormal{ann}(\hat y)&\hat 0\end{array}\right]=\left[\begin{array}{cc}\hat 0&\hat 0\\\hat0&\hat 0\end{array}\right],\]
because $\textnormal{ann}(\hat y)=\hat 0$. Therefore the righthand side of (\ref{eq71}) is equal to
\[\left\{\left.\left[\begin{array}{cc}\hat h_1&\hat y\hat
h_2\\\hat x\hat h_2&\hat h_1-\hat y\hat
h_2\end{array}\right]\right|\ \hat h_1,\hat h_2\in
F[x,y]/\langle x^2\rangle \right\},\]
which does not contain the matrix $\left[\begin{array}{cc}\hat x&\hat x\\\hat 0&\hat 0\end{array}\right]$. However, direct verification shows that
\[\left[\begin{array}{cc}\hat x&\hat x\\\hat 0&\hat 0\end{array}\right]\in \textnormal{Cen}(\widehat B).\]}
\end{Example20}
\vskip 0.3 cm
In the following example we will see that for every $n\ge 3$ and for any UFD $R$ and ideal $I$ such that $R/I$ is a ring with zero divisors, there is a matrix $B\in M_n(R)$ for which we do not have equality in Proposition \ref{Lemma5}.

\newtheorem{Example21b}[Theorem2.1]{\bf Example}
\begin{Example21b}\label{Example21b}
\textnormal{Let $R$ be a UFD and let $I$ be an ideal in $R$ such that $R/I$ has zero divisors. Thus suppose that $\hat d\hat d'\in R/I$, $\hat d,\hat d'\neq\hat 0$ and $\hat d\hat d'=\hat 0$. Now let $B=\left[\begin{array}{ccc}0&d&1\\0&0&1\\0&0&0\end{array}\right]\in M_3(R).$ Note that $d\neq 0$ since $\hat d\neq\hat 0$. Because the characteristic polynomial of $B$ is equal to the minimum polynomial of $B$ it follows from Theorem~\ref{Theorem2.1} and Lemma \ref{newlemma} that
$\textnormal{Cen}_{M_3(R)}(B)=$\[\left\{\left.a\left[\begin{array}{ccc}0&0&d\\0&0&0\\0&0&0\end{array}\right]+b\left[\begin{array}{ccc}0&d&1\\0&0&1\\0&0&0\end{array}\right]+c\left[\begin{array}{ccc}1&0&0\\0&1&0\\0&0&1\end{array}\right]\right|\begin{array}{l}a,b,c \textrm{ are elements}\\\textrm{ of the quotient}\\\textrm{ field of $R$.}\end{array}\right\}\cap M_3(R),\]
and so every matrix in $\Theta(\textnormal{Cen}(B))$ has $\hat 0$ in position $(2,1)$. Furthermore, using the notation in Proposition \ref{Lemma5} we have
\[[\mathcal A_{ij}]=\left[\begin{array}{ccc}\hat 0&\hat 0&R/I\\\hat 0&\hat 0&\textnormal{ann}(\hat d)\\\hat 0&\hat 0&\hat 0\end{array}\right].\] Hence, every matrix in $\Theta(\textnormal{Cen}(B))+[\mathcal A_{ij}]$ has $\hat 0$ in position $(2,1)$. However, direct multiplication shows that\[\left[\begin{array}{ccc}\hat d'&\hat 0&\hat 0\\\hat d'&\hat 0&\hat 0\\\hat 0&\hat 0&\hat d'\end{array}\right]\in \textnormal{Cen}(\widehat B),\] and so equality in Proposition \ref{Lemma5} does not hold in this case.
Now, again let $R$ be a UFD and let $I$ be an ideal in $R$ such that $R/I$ has zero divisors. Let us consider the matrix \[B'=\left[\begin{array}{c|c}\begin{array}{ccc}0&d&1\\0&0&1\\0&0&0\end{array}&\bigcirc\\\hline\bigcirc&\bigcirc\end{array}\right]\in M_n(R).\]Then
\[\textnormal{Cen}(B')\subseteq\left[\begin{array}{c|c}\textnormal{Cen}(B)&R/I\\\hline R/I&R/I\end{array}\right]\qquad
\textrm{and}\qquad[\mathcal A_{ij}]\subseteq\left[\begin{array}{c|c}\begin{array}{ccc}\hat 0&\hat 0&R/I\\\hat 0&\hat 0&\textnormal{ann}(\hat d)\\\hat 0&\hat 0&\hat 0\end{array}&R/I\\\hline R/I&R/I\end{array}\right].\]
Since \[\widehat A:=\left[\begin{array}{c|c}\begin{array}{ccc}\hat d'&\hat 0&\hat 0\\\hat d'&\hat 0&\hat 0\\\hat 0&\hat 0&\hat d'\end{array}&\widehat\bigcirc\\\hline \widehat\bigcirc&\widehat\bigcirc\end{array}\right]\in \textnormal{Cen}(\widehat{B'}),\]but clearly $\widehat A\not\in\Theta(\textnormal{Cen}(B'))
+[\mathcal A_{ij}]$, equality in Proposition \ref{Lemma5}, for these cases, does not hold.
}
\end{Example21b}

\bigskip

\end{document}